\documentclass[12pt]{article}
\usepackage[latin1]{inputenc}
\usepackage{amsmath}
\usepackage{amsfonts}
\usepackage{amssymb}
\usepackage[usenames,dvipsnames]{color}%\usepackage[dvips]{color}
\usepackage{graphicx}
\title{\bf  The UK financial mathematics M.Sc.}
\author{ 
{\sc A.E. Kyprianou
}
\\
Department of Mathematical Sciences, University of Bath.\\
\bigskip
}

\begin{document}

\maketitle

%\begin{appendix}
\begin{center}{\bf Abstract} \\
%\texttt{Not for Circulation.}
\end{center}

\bigskip

\noindent  Postgraduate taught degrees in financial mathematics have been booming in popularity in the UK for the last 20 years. The fees for these courses are considerably higher than other comparable masters-level courses. Why? Vendors  stipulate that they offer high-demand, high-level  vocational training for future employees of the financial services industry, delivered by academics with an internationally recognised  research reputation at world-class universities. 

\medskip

\noindent  We argue here that, as the UK higher education system  moves towards a more commercial environment, 
 the widespread availability of the M.Sc. in financial mathematics exemplifies a  practice of following market demand for the sake of income, without due consideration for the broader consequences. Indeed, we claim that, as excellent as  such courses can be in intellectual content and delivery, they are mismatching needs and expectations for such education and confusing the true value of what is taught.
 
  \medskip
 
 \noindent
 The story of the Mathematical Finance MSc  serves as a serious case study, highlighting some of the  incongruities and future dangers of free-market education.

%this practice exemplifies an encroaching sense of  opportunism  within academia that is contributing to a bubble effect. 
%\end{appendix}

\section*{The UK financial mathematics M.Sc.}

The last three decades have seen an unprecedented explosion of activity in a new sub-discipline of mathematics: {\it financial mathematics}. The emergence of this field of research and study  has paralleled the manifold expansion of the financial services industry and, accordingly, has transformed  the role that mathematics departments now play in feeding the financial services industry. % with highly educated quantitative analysts.

\medskip

\noindent For many university mathematics departments, the volume of demand for taught courses that have been set up in the direction of financial mathematics, in particular at the postgraduate level, has proved to be a gift from the gods. In the UK, the country to which we henceforth restrict the majority of our forthcoming discussion, the {\it Master of Science (M.Sc.) in financial mathematics} has been pioneered for the best part of the last twenty years. Such masters programmes are marketed as premium postgraduate education, with the allure of a well-paid career in the banking industry thereafter. Moreover, they are often populated by large numbers of foreign students unable to access this type of specialised education in their country of origin. Accordingly, these study programmes are generally priced  with some of the highest annual  tuition fees in the UK. (To give a little perspective, quite a few universities are now charging well in excess of what one would pay for a top-end, brand new family-sized car.)  For some universities, the fees are found to be  double, treble or even quadruple  those  of many other masters programmes on the same campus, in the same faculty or even in the same department. A brief glance at the table in the Appendix A reveals some eye-watering statistics. Appendix B illustrates further the severity of the differential between the cost of undertaking a regular  M.Sc. in Mathematics (or a related course with a very strong mathematical component) compared to an M.Sc. in Financial Mathematics. For home/EU students, the situation is considerably worse than for overseas students who, traditionally, have always paid significantly higher fees to attend UK universities. Indeed one might even question from the data presented whether home/EU students are being priced out as there is a deliberate focus on the overseas market where more money is to be made. Nonetheless, even for overseas students, there is a clear and significant pricing discrepancy with `nearest-neighbour' mathematics masters programmes.

\medskip

%\noindent %Unlike the issues surrounding undergraduate fees in the UK, which, as of

%In  2012, undergraduate fees in the UK were dramatically increased threefold for home and EU students at almost all of the leading universities. Nonetheless, with the exception of a tiny number of private universities, they remain capped at a maximum of \pounds 9,000.  In contrast,   

\noindent For many years, the setting of fees for postgraduate taught courses has been  left to the discretion of the awarding institution. It is assumed  that the principle of the free market should come into play. Universities are expected to set fees that represent a balanced view on reputation, expertise, running costs and demand. As British academia moves into an age where the commercialisation of education is becoming ever more the norm, we examine whether things have gone too far in the case of the M.Sc. in financial mathematics, thereby highlighting some of the  incongruities and future dangers of free-market education in the UK.

%But to what extent do universities, consciously or unconsciously,  exaggerate the value of the education they offer in this respect?  There is growing debate and concern that a bubble-like effect is emerging in this `education market'. The   M.Sc. in financial mathematics provides an interesting and quite telling case in point. 

\section*{Niche mathematics}

\noindent In 1997, Robert Merton and Myron Scholes were awarded the Nobel  Prize for Economics for their work on the use of probability theory and its connections to partial differential equations, in order to show how it was possible to explicitly  hedge  and hence price  a whole family of  option contracts in an idealised market. One of the key innovations in their work was the   translation  of rational behaviour of market agents into a rigorous mathematical framework.
Although their model of a financial market was an idealised version of reality, their calculations involved a  remarkable and stunning  application of one of the most important developments in the theory of probability:  {\it stochastic calculus}. This was a truly exciting development that provided a rare example of where a highly non-trivial mathematical theory could be used  to address  a  real-world application. Moreover, it  offered  the possibility of much more sophisticated computations and, eventually, opened  up the sub-discipline that we now call {\it financial mathematics}. 

\medskip

\noindent This Nobel prize-winning work dated back to the mid-seventies, and pertained  to an original publication that Scholes had co-authored with Fisher Black in 1973. (Black died in 1995 and could not be awarded the prize  posthumously.) Black, Scholes and Merton were all financial economists and, whilst Black and Scholes had some initial difficulties publishing their work, momentum in the development of the so-called Black-Scholes-Merton theory developed through the late seventies. By the early 1990s this theory had captured the attention of the mathematical mainstream. Thanks to a general sense of awareness, in part brought about by  the award of the Nobel Prize, a much broader perspective on the original contributions of the financial economist trio had been realised by the late 1990s and the amount of research being published in mathematical journals and books  blossomed dramatically. By the turn of the Millennium,   financial mathematics was well and truly established as a field of research. Today, it claims the right to be recognised as a mainstream sub-discipline of mathematics and is supported by a huge international community of researchers.

\medskip

\noindent In parallel to (as well as  stimulating)  these academic  developments, the global derivatives markets had significantly swollen in size thanks to a process of deregulation that began over 30 years ago, predominantly in the US and UK.  By the 1980s, traders were  selling in large volume complex contracts `over the counter', which were generally known as {\it financial derivatives}. For an up-front premium,  a financial derivative would allow investors to collect potential future financial returns, or perhaps no return at all, depending on how future markets played out. Option contracts, such as those originally studied by Black, Scholes and Merton, were relatively straightforward examples of such derivatives. By the mid 1990s, for good or for bad, the volume and variety of financial derivatives that were being traded reflected the degree of sophistication with  which traders were using these products to build complex investment portfolios. 
As markets are unpredictable, each financial derivative carries an individual risk and, accordingly, portfolios of derivatives carry compounded and interactive risks. With such large-scale and complicated trading taking place, financial institutions sensed the increasing need for in-house mathematicians to quantify the risks involved. %This meant being  able to value  the current and projected value of investment portfolios. 
And so, the number of positions for so-called {\it quantitative analysts}, or {\it quants} for short, along with other supporting roles of a numerate nature, has exploded in number over the last 20 years or more. 

\medskip

\noindent Together with the US, UK academia was quick to respond to this rapidly expanding job market and to identify a rather unique opportunity to set up specialist postgraduate mathematics teaching programmes that, in some sense, also served as high-level vocational training for the quantitative financial services industry. Here was a self-reinforcing niche market. As more financial mathematics graduates  proliferated into various quantitative support roles in the financial sector, the more the next generation of students perceived the necessity to engage with such training. The M.Sc. in financial mathematics has become the iconic intermediate goal of many young women and men, particularly those with  a first degree in mathematics, statistics  or physics, whose long-term aspiration is to  `work in banking'.

\section*{Cash cow}

\noindent The very first UK-based M.Sc. in financial mathematics made its debut at the University of Warwick in 1996. A similar programme was launched the following year jointly by the University of  Edinburgh and Heriot--Watt University in the Scottish capital. According to their webpages, both programmes still enrol  very healthy numbers today: approximately 50 in Warwick every year and 25-35 in Edinburgh every year. This is even more impressive when one takes account of the number of competing courses that have sprung up in the UK  in the meantime, not to mention in other countries around the world. As can be seen from Appendix A, there are currently over 30 M.Sc. courses on the topic of financial mathematics,  that can be found the length and breadth of the UK. The majority of annual fees for these masters programmes exceed the upper cap on undergraduate fees (\pounds 9,000). Around two thirds of them  exceed \pounds 15,000, going as high as  \pounds 31,000, for home and EU students. Whilst it is a long-standing practice  in the UK higher education system that a higher rate is charged  for non-EU students, the majority of masters in financial mathematics do not discriminate between whether a student is from the EU or not. A back-of-the-envelope calculation shows that a department with a well maintained M.Sc. degree of this kind should easily manage to bring  in well over half a million pounds a year from it. Some programmes currently in existence are more than capable of bringing in well over a million pounds per year. %These figures are even easier to achieve when the majority of the cohort are non-EU students.

\medskip

\noindent The lure of this kind of income for mathematics (and partner) departments has proved to be far too attractive to ignore in light of the  financial strain that British academia is under. Over the last 20 years, encouraged by the  government, UK universities have experienced a massive programme of expansion, opening their doors to many more home and international  students.
Naturally, state funding has  spread thinner and become seemingly less generous. Student grants have been replaced by student loans,  study fees for undergraduate home students have been passed from the state to the student,  national research funding council scholarships for taught postgraduate courses have all but disappeared and research grant capture has become increasingly more competitive. Universities have found themselves in a position where, in order to survive  in the modern academic climate, it has become  necessary to explore a whole variety of new sources of income.

\medskip

\noindent Laying on low-cost, easy-to-populate taught postgraduate degrees with premium rate fees is a straightforward way of generating income. Pursuing the ever-expanding, higher fee-paying international student market is another. 
Running an M.Sc. in financial mathematics has proved to be very successful  in both respects for many universities.

\section*{Value for money}

% \noindent  On the back of this kind of postgraduate teaching income, departments have been able to justify the creation of new academic posts to support teaching in financial mathematics. Such appointments  are also used to bolster the perceived research profile in the field of financial mathematics. %, thereby giving credibility to the the associated financial mathematics masters programme. 
 %An increase in research and teaching activity  in financial mathematics has brought about the inevitability that some students have undertaken  research leading to a Ph.D. at the same or another university. This natural progression  is actively encouraged and forms an important part of the process of  maintaining credibility when marketing postgraduate education in financial mathematics. It is often sold as a  mark of quality that significant and internationally visible  research activity can be  seen behind the teaching machinery.

%\medskip

%\noindent But, realistically speaking, how well are academics in touch with the very industry that they claim to be providing specialised postgraduate training for?

\medskip

%\noindent Ask any academic involved in the teaching of postgraduate financial mathematics if there is any interaction with industry and most will have something constructive to say. Some may refer to former students with whom they still hold dialogue, some may refer to joint seminars with industrial participation, some may even point towards funding that has streamed into the department as a consequence of industrial interaction. Certainly, at one extreme, Oxford University can boast the existence of the Oxford Man Institute, an interdisciplinary research centre which is heavily subsidised  with hard cash from industry. On the other hand, nationally speaking, it is likely that only a few of those involved in running and teaching on  postgraduate masters programmes in financial mathematics have regular interaction with the banking industry, or indeed awareness of the kind of tasks, routines and challenges that are undertaken by quantitative analysts on a daily basis.
%\medskip

\noindent  How is value for money of premium-fee M.Sc. programmes justified? {\it Specialism}, {\it opportunity} and {\it demand}  are three of the words that many offering such education cite.

\medskip

\noindent For the large majority of financial mathematics masters programmes on offer, it is not unusual that a proportion of what is taught is standard material which can be found  in third or fourth year undergraduate programmes. Indeed, it is not uncommon that, in the interests of the economy of scale,   masters students will find themselves sharing the lecture theatre with undergraduates and postgraduates of other masters programmes for some of their modules.  %Other modules may be exclusively taught to the M.Sc. students often illuminated in an appropriate light  to accommodate for the particular fee-paying audience. For example, it is not uncommon to teach a course on stochastic calculus whose punchline at the end is an exposition of its application in the original Black-Scholes-Merton Nobel-prize-winning theory.
% One similarly  sees other courses which are hosted by the local management or business school which pertain more to quantitative economic theory. 
Specialist material, exclusive to the masters programme itself, is also to be found. 
For example, it is normal in the UK  that the successful completion of an  M.Sc. requires the submission and satisfactory assessment of a  summer dissertation. Here the student has the opportunity, over a period of around 3-4 months, to be guided through an extensive programme of reading and research, addressing considerably more advanced theory and  application. However, most of what is offered up in this respect is considered `specialist' because it addresses current trends in academic research rather than what is actually used or needed in the industrial roles that  the  majority of these masters  graduates go on to occupy.  This discrepancy is an important point that is worth emphasising.

\medskip 

\noindent %The volume of published research  in financial mathematics is vast. 
A surprising fact concerning modern research in financial mathematics  is that, although presenting substantial technical challenges,  it does not  find its way into industry as often as one might expect. Many of the underlying models which are studied in theory are simply inappropriate or impractical to work with. It is not uncommon for the academic financial mathematician to work from the safety of an axiomatic treatment of a financial scenario in order to investigate phenomenological behaviour or to  develop a mathematical method within that context.  The discrepancy with reality comes about because, unlike for example physical systems,  it is difficult to reasonably model  the mechanisms that drive the randomness in generic financial markets by appealing to fundamental universal laws. %The axiomatic approach taken in the vast majority of academic research, as complicated as it may be, is to simple to have significant predictive integrity in real life. %Many who  associate their research with financial mathematics are concerned with  phenomenological properties of simplified models 
%For example, just the influence on market behaviour of investors' psychology alone  is something that in itself no one has yet managed to find a reasonable mathematical treatment of.

\medskip

\noindent It is not true, however, that academic research in financial mathematics has no impact at all. There are some very important successful examples attributed to quite exceptional mathematicians. These individuals  have carved out relatively rare careers, interacting closely and consistently with banks and   influencing  in-house banking activity, sometimes from within an academic position, sometimes as an employee of a bank and sometimes both. But such individuals are few and far between and the level of mathematical experience involved goes far beyond what can be established on the back of an M.Sc. in financial mathematics.
Moreover, it is usual that their interaction takes place with some of the bigger banks or specialist consultancies that have a unique inner core of `super-quants' who are educated and experienced as academic researchers well beyond the level of a Ph.D. % (usually in possession of a Ph.D. and even with some postdoctoral academic research experience). Indeed, it would not unusual for such super-quants  to publish regularly in academic journals. 
%
%
%\medskip 
%
%
%
%\noindent 
%It is therefore  `specialism' on offer in many of these postgraduate programmes 
%In general, interaction between academics and the community of quantitative analysts is quite poor. 
%Nonetheless, this not the typical picture of how academics and quants interact with one another, not by far. 

\medskip

\noindent %Most of what academics enjoy associating themselves with by way of `specialism' on the theme of financial mathematics never gets even remotely close to the average quant's desk. 
As far as masters-level vocational training is concerned, the reality for most  academics is that  many are not up to date on the kind of challenges and routines that are undertaken en masse by industrial quantitative analysts on a daily basis.  Moreover, few can claim to have recently seen the inside of a bank and interacted in a mathematical capacity with its quants. %, even for the large majority of academics who actively publish in the field of financial mathematics.
  Therefore, one has to question the extent to which  the academic specialism that is injected into postgraduate financial mathematics taught courses should  be valued {\it over and above} any other numerate discipline  taught at this level, in terms of vocational training. %, just because of a collective association with a profession offering high salaries. 

\medskip 

\noindent That said, there is still the issue of what employers want. Webpages, posters and brochures for postgraduate education in financial mathematics often advocate the value of the course in question as useful in preparing for a job in the financial industry and point towards some of the many employers that their graduates have gone on to work for. But, realistically, to what extent are banks interested in hiring individuals with such specific masters qualifications?   It would be more appropriate if marketing for these programmes could use hard data to indicate the true demand for such education {\it from within the industry itself}. In the experience of the author  of this article, an academic who {\it has} spent a modest amount of time inside a number of banks and financial institutions engaging with their quantitative analysts, very few banks seem to care whether or not a student is in possession of an M.Sc. in financial mathematics.

\medskip

\noindent Indeed, hiring strategies for quants (the industrial job whose technical requirements most closely match the training on the typical Financial Mathematics M.Sc.)  are  based around looking for sharp  problem solvers who can deal with the fast pace of work and learn and develop in-house technical procedure. The familiarity with certain concepts and academic specialisms that an M.Sc. in financial mathematics offers, some of which are not directly relevant to working practice, is not necessarily deemed to be of huge value to the quant employer, who would expect that the smart hire  can  pick  up what is needed anyway once in post.
Often, but not always, this means that it is those who are trained to the level of a Ph.D. in a numerate discipline (not necessarily financial mathematics) that have a competitive advantage as far as quant jobs are concerned. If nothing else, this is because of a self-selecting bias through their ability to undertake a Ph.D. in the first place. 

\medskip
%\noindent Unlike, for example, the world of actuaries or professional statisticians, where a specialist masters-level education does bear direct relevance to the vocational training of the future employee, 

%Once again, we claim there is a case to be made here that   a masters education in financial mathematics does not necessarily   make a difference in terms of employability {\it over and above} any other undergraduate or postgraduate numerate degree; that is, at least for those students who do not wish to pursue their education any further thereafter. 

%\medskip

\noindent This does not contradict the observation made by vendors of financial mathematics postgraduate education that many of their graduates do go on to work in the finance industry. It just means that the technical needs of the roles they are hired into (for example careers in risk management or trading, which have relatively soft graduate technical numerical requirements) are not specific to an M.Sc. in financial mathematics. Moreover, self-selection again plays a role here; commitment to an expensive one-year postgraduate certificate is a strong affirmation of the desire to end up working in the financial services industry.

\medskip

\noindent Another justification for the value for money that can be given, but only for a handful of M.Sc. programmes in financial mathematics (see Appendix A), is that they offer the opportunity for   a summer placement in industry or a summer research project in collaboration with industry. However, with the exception of one university, it is not entirely clear from publicly available on-line course prospectuses exactly how many students on each programme have this experience accessible to them during the summer months. Certainly, if mentioned at all, the choice of language one tends to find regarding the exact nature of such opportunities within the degree is rather vague and non-committal.  

\medskip

\noindent More generally, participating in an M.Sc. in financial mathematics will almost certainly give students the opportunity to enter job interviews in the financial services industry at large (not just quant jobs) with more confidence. But it seems somewhat anecdotal to claim that this would {\it drastically} improve one's chances of employment over and above any other masters level education in mathematics given that the core theory that supports such masters programmes is relatively easy to  access in most mathematics undergraduate and taught postgraduate programmes anyway. Certainly there is no apparent data (at least not supplied by vendors of financial mathematics masters programmes) to support such a claim. In particular, the kind of data that is desired here would certainly need to take account  of the phenomenon of self-selection as noted above.

\medskip

\noindent All of the above arguments aside, one still cannot discount the fact that there is ongoing demand for these masters courses from the students themselves, not to mention many more financial mathematics bachelor programmes that are also in existence (priced in line with standard undergraduate fees). More broadly, there is at least one commercially operated programme of courses which could be argued to be in competition with university education. A Certificate in Quantitative Finance is provided by Fitch Learning in collaboration with the Wilmott Forum, which currently charges \pounds 11,950 for a six-month programme with six specialist modules. 

\medskip

\noindent As long as students are keen to pay the current costs involved in postgraduate training in financial mathematics, then one could assert that the  product speaks for itself. But, realistically speaking, is this an acceptable justification for the extent of the pricing differential observed in universities? 

\medskip

\medskip

\noindent Amidst the  prevailing culture shift  in UK  academia to a more corporate approach, one has to ask whether `opportunistic education' has started to appear. The packaging, marketing and selling of education is part and parcel of  the many challenges that universities are now confronted with as a matter of survival. Against this backdrop, one has to question the motivation behind the scramble  in recent years by institutions wanting to offer  postgraduate education in financial mathematics whilst imposing relatively (in some cases, exceptionally) high fees, when, given the arguments above, there isn't  necessarily clear evidence of  value for money at that tariff.  There are grounds to feel concern  that there is an element of implicit cynicism around the running of such courses. Said another way, there is a danger that the answer to the question {\it `why are fees for financial mathematics postgraduate taught courses so  high in relative terms?'} may simply be: {\it `because universities can get away with it'.}

\medskip

\noindent Would any other public service industry tolerate or be tolerated for such an approach to its pricing policy?

\section*{Responsible  education}

Financial mathematics is beautiful. It has proved itself again and again, showing amazing robustness in the way that financial and economic concepts and problems can be reformulated in an abstract context and addressed using a wide variety of different mathematical tools. Moreover, there can be no doubt that heightened interest in financial mathematics has helped to motivate  some of the major waves of development in probability theory over the last 30 years.
 Nonetheless, financial mathematics has its critics. There are those who feel that the intelligentsia should not promote a false sense of security in the use of mathematics where there  are known to be inaccuracies, especially  within an economic system that has led to an extreme, and ultimately devastating, financial breakdown.

\medskip

\noindent But financial mathematics is not about making money and it is not about promoting capitalism (or any other socio-economic ideal for that matter). Financial mathematics is a science, it is a theory that aims to quantify risk and randomness within economic scenarios, with a view to characterising {\it rational} behaviour, thereby promoting {\it fairness}. In the experience of this author, most mathematicians who advocate the use of this theory (both within academia and within the financial services industry) have a sincere interest in seeing it used to improve the quality of decision making within financial markets, leading to a more accountable and responsible system. %The politics of how those markets are operated and whether they should even exist is a much larger debate. %The inclusion or exclusion of mathematics in that debate is of minor influence.  % and ultimately aspire to being a part of a process of quantifying regulation of and restraint from the kind of wild risk-taking that has occurred in the past.

\medskip

\noindent In this respect, financial mathematics {\it is} an important intellectual discipline and deserves to be taught at depth in an academic environment.  There is a lot to share and a lot to think about. And like all other scientific disciplines, it deserves to be delivered in a responsible way. In this respect most of the masters programmes in existence are extremely rewarding in their intellectual content. Indeed, some of the programmes listed in Appendix A are  exemplary  in their academic delivery, for example, in  turning their attention to the mathematics of post-crash finance. But, whilst this is no doubt relevant to the issue of value for money, there is a much bigger picture to consider here.

\medskip

\noindent %Given the remit of the `financial mathematician' as one who seeks rational behaviour, it seems overly ironic that universities %(whose interest in setting premium fees on such masters programmes has partly come about due to  loss of income from state funding following the financial crisis)  have allowed the rapid inflation of postgraduate fees.
There is a real danger that the M.Sc. in financial mathematics has  become an iconic symbol within UK academia of an overzealous attitude towards  the commercialisation of education.
Prospective students are ill informed about the relationship between research in financial mathematics, which guides what is taught on  masters programmes devoted to this theme, and what is actually needed in the banking  industry for the different types of numerate jobs. Marketing for masters programmes in financial mathematics often confuses what is needed by an elite core of quantitative analysts in large banks and consultancies (who typically look far beyond the training that a UK M.Sc. in financial mathematics has to offer) with other softer quantitative roles in banking (for which a whole array of other quantitative undergraduate and postgraduate qualifications are equally relevant).  

\medskip

\noindent The fee structure  of the UK financial mathematics M.Sc.  sets an undesirable precedent for other degrees in  the future. Take, for example, postgraduate taught degrees in statistics. Unlike the case of an M.Sc. in financial mathematics, one can make a much stronger case that there {\it is} directly relevant high-level vocational training taking place in a statistics M.Sc. Moreover, such a qualification genuinely {\it is} considered to be needed by a whole variety of industries that {\it do} hire masters-level students into well-paid jobs as statisticians. What is now stopping universities charging significantly higher fees than the average for these degrees? Indeed,  there is plenty of demand in this field too.
%
%\medskip
%
%\noindent 
 Do we want our academic system in the UK to fall into the  trap in which academic excellence is confused with costly branding? As a quick google will confirm, there is already much debate  of an emerging bubble in the education market, in particular in the US.  If this is the case, would it not be a sad irony if {\it financial mathematics} found itself at the head of such a market in the UK?
% Once again, the fact that it is financial mathematics on the front line here is far too much of a painful irony.

\medskip

\noindent %Is it time for universities to take tighter control of the marketing of these programmes and consider more carefully the effects of current pricing regimes?

\medskip

\noindent Arguing the case of  a sellers' market for this kind of education seems the most dangerous justification of all for premium fees. Given the significant differential  in fees with nearest-neighbour mathematics postgraduate degrees for home-EU students  at a large number of institutions and given the large proportion of non-EU students that are enrolled every year, universities could potentially risk exposing themselves to an accusation of `predatory education' directed at wealthy foreign students. If the principle of appealing to `demand' to justify the setting of exceptionally high fees becomes a more common practice, would the system not attract  the uncomfortable prospect of external regulation?

\medskip

\noindent For the specific case of the M.Sc. in financial mathematics, justifying premium fees on the basis of there being a sellers' market also opens up another genuine concern. On completion of their studies, some overseas students are returning to developing and less regulated economies where it is unclear to what extent a qualification of this kind, marketed and sold to students the way that it is, is being put to use. This aspect alone may prove to be the most frightening consequence of all in  the  financial mathematics postgraduate education frenzy. Who is checking these things?

\medskip

\noindent Twenty years ago, when the profession of the quantitative analyst was still relatively young, the M.Sc. in financial mathematics meant something completely different to what it does today. The financial services industry has matured from where it was and, at least in the main financial hubs  (London, New York, Frankfurt, Paris, Hong Kong, Singapore, Tokyo, etc.), it has a much clearer sense of identity and a clearer understanding of how it is prepared to interact with high-end mathematics.
The role that universities can play in vocational training for this profession has accordingly changed.
As academics, we need to take a closer look at  the data available to correctly understand, quantify and qualify the true value of what is currently on offer by way of taught postgraduate education in financial mathematics. %Moreover, we need to understand the effect we are having outside of academia as a consequence of the theories we teach to our graduates. 

%\medskip

%\noindent %The story of masters education in financial mathematics is but a snapshot of a much bigger picture. The selling of education has become an industry in itself. Is what we see for financial mathematics the shape of things to come?
% and, like all industries some degree of regulation is needed. Historically, the examination of all degree programmes, both undergraduate and postgraduate are monitored by an external examiner, typically an academic from another university whose individual tenure in the role of observer lasts no more than three years.    But external examiners are only invited to comment on the quality of academic content and the process and delivery  of examinations. Decisions such as determining the value for money are made at the university's own discretion. Who is there to protect the interests of the students when it comes to the setting of fees? 
%Do we want our academic system in the UK to fall into the  trap in which academic excellence is confused with costly branding?

\medskip

\noindent   Conversely, {\it if} UK university education in this field is to be subjected to the full force of the  free market ideal, then is it clear that  the students, i.e. the consumers, are correctly informed when choosing to enroll for a M.Sc. in financial mathematics? Students should demand greater justification when premium fees are involved. 
How do these fees break down?  Can the university really quantify the value it will have in their future career on demand (other than a fancy poster, brochure or webpage with quotes from satisfied students)? Can a university's claims be independently verified? What is the relative additional value for money over and above nearest-neighbour masters programmes whose graduates also feed the financial services industry, and whose fees already take  account of the fact that graduate education has value to employers at large? %To put things in perspective, some universities are now charging for an M.Sc. in mathematical finance well in excess of what one would pay for a top-end, brand new family-sized car. Considering the amount of thought and consideration of facts and figures  (which are readily supplied by car manufacturers and their associated vendors) one would  put into investing in the latter, why would the situation be any di
%
%\medskip
%\noindent Rather than capable students  being encouraged to fork out large amounts of money on the basis of a university's  self-perceived worth, would it not be better that {\it universities} are priced out of chasing cash cows by the students themselves? This would require prospective students to take  a more aggressive approach in the educational market place, placing higher demands on the justification of cost as well as directly haggling the price down by playing one vendor off against another.  

\medskip

\noindent
After all, 
a fundamental concept that  students will learn on an M.Sc. in financial mathematics is that rational pricing forbids the occurrence of arbitrage.

\section*{Acknowledgements}

\noindent The motivation for this document comes from  
more than 15 years of the author's experience in teaching, managing and externally examining M.Sc. programmes in financial mathematics  across  a wide spectrum of universities, both nationally and internationally. %He has also managed  one of the very first UK-based M.Sc. programmes at Edinburgh University in the late 1990s. 
Parts of this document are  inspired by many discussions with colleagues over the aforesaid period in light of the growing commercialisation of education in the UK. More recently, whilst writing this article, the author has had a number of extensive, candid and extremely helpful  discussions with stakeholders  in postgraduate financial mathematics education, in both academia and industry.  I am most grateful for their input.

\newpage

\section*{Appendix A}
{\footnotesize The following table lists current postgraduate taught courses in the UK with emphasis on financial mathematics, ranked in decreasing order of fees for home/EU students. This data has been collected from on-line, publicly available  information provided by each of the institutions. Fees indicated are for the academic year 2014-2015 or 2013-2014 where information for the former is not yet available. Figures, expressed in k\pounds, have been rounded to one decimal place. The list of postgraduate courses that relate to financial mathematics is  more extensive than the one presented below, but we have chosen to include only the ones which present significant mathematical content. %. The judgement to exclude a course from this list is based  on the volume of its mathematical content.} 
The fourth column, `P', indicates whether there is mention   on the programme webpages of the possible availability of  of an industrial placement or industrially oriented research project during the summer dissertation period. In this column, `N' means there is no mention, `V' means there is vague mention, `E' means there is explicit mention of opportunities (without quantifying how many students will have access to this opportunity or whether there are any guarantees) and `Y' means that a strong commitment to engaging the majority of students with an industrial placement or industrially oriented research project has been indicated.

\[
\hspace{-0.4cm}\begin{footnotesize}
\begin{array}{| l | c | c |c|l|}
\hline
 \text{
University} & \text{  home/EU (k\pounds)} & \text{ non-EU(k\pounds)} & \text{P} & \text{M.Sc. Course name}\\
\hline
\hline
\text{Warwick}   & 31&31 &\text{E}& \text{Financial Mathematics}\\
\hline
 \text{Oxford}    & 27.5 & 31.1 &\text{N}^{1}& \text{Mathematical Finance (part time)}
\\
\hline
 \text{Imperial}  & 27 & 27 &\text{Y}& \text{Mathematics and Finance}
\\
\hline
\text{Oxford}    & 24.9 & 24.9 &\text{N}&  \text{Mathematical and Computational Finance
}
\\
\hline
 \text{City}  & 23 & 23 &\text{V}&  \text{Financial Mathematics}
\\
\hline
 \text{LSE}  & 23 & 23 &\text{N}&  \text{Financial Mathematics}
\\
\hline
\text{LSE} & 23 & 23 &\text{N}&  \text{Risk and Stochastics}
\\
%\hline
%8 &\text{LSE}   & 23 & 23 & \text{Statistics (Financial Statistics)}
%\\
\hline
 \text{King's}   & 22.7 &  22.7 &\text{N}& \text{Financial Mathematics}
\\
\hline
 \text{Manchester}   & 21.4 & 21.4 &\text{E}&  \text{Quantitative Finance: Financial Engineering }
\\
\hline
 \text{Reading}   &19.8  &  19.8&\text{N}& \text{Financial Engineering}
\\
\hline
 \text{UCL}   &19.6 &21.7 &\text{N}& \text{Financial Mathematics}
\\
\hline
 \text{Birkbeck}   & 19 & 20 &\text{N}&  \text{Financial Engineering}
\\
\hline
 \text{Edinburgh/Heriot--Watt}    &19 &22 &\text{E}&  \text{Financial Mathematics}\\
%\hline

%\text{Heriot--Watt}  & 18.3 & 18.3 & \text{Quantitative Financial Engineering}\\

\hline
\text{Brunel}   & 17.5 & 17.5 &\text{V}&  \text{Financial Mathematics}
\\
\hline
 \text{York}   & 17.4&22.2 &\text{N}&  \text{Financial Engineering}
\\
\hline
\text{York}   & 16.5&21 &\text{N}&  \text{Financial Mathematics}
\\
\hline
 \text{Queen Mary}   &15.5  &  18&\text{N}&  \text{Mathematical Finance}
\\

\hline

 \text{Birmingham}  & 15.5 & 15.5 &\text{N}&  \text{Mathematical Finance}
\\

\hline

 \text{Liverpool}   & 15.4 & 15.4 &\text{N}&  \text{Financial Mathematics}
\\

%\hline
%\text{Lancaster} &  15&19   & \text{Quantitative Finance}
%\\

\hline

 \text{Exeter}    &11.9 &19.5 &\text{N}&  \text{Financial Mathematics}\\
\hline
 \text{Manchester}  & 10.7 & 18.4 &\text{V}&  \text{Mathematical Finance}
\\

%\hline
%\text{Cambridge} & 9.6 &18.7  & \text{Financial Engineering (M.Phil.)}
%\\
 \hline
\text{Leeds}   & 9 & 19 &\text{N}&  \text{Financial Mathematics}
\\

\hline
 \text{Leicester}   &8.5 &16 &\text{E}&  \text{Financial Mathematics and Computation}
\\

\hline
\text{Glasgow}  & 8.3& 17.2 &\text{N}&  \text{Financial Modelling}\\
%\hline
% \text{Glasgow}   & 8.3& 17.3 & \text{Quantitative Finance}\\

\hline
 \text{Strathcylde}   &  7&13  &\text{N}&  \text{Quantitative Finance}
\\

\hline
 \text{Loughborough}   & 6.3& 13.8&\text{N}&  \text{Mathematical Finance}
\\

\hline
 \text{Birmingham}   & 6 & 13.7 &\text{N}&  \text{Financial Engineering}
\\

\hline
 \text{Sheffield}   &  6 &15.4  &\text{N}&  \text{Statistics with Financial Mathematics}
\\

\hline
 \text{Nottingham}   & 6& 13.7 &\text{N}&  \text{Numerical Techniques for Finance}
\\

\hline
 \text{Swansea}   & 5&12.5 &\text{N}&  \text{Mathematics and Computing for Finance}\\

\hline

\end{array}
\end{footnotesize}\]
$^1$ This programme is designed for individuals who are already employed in the financial services industry.

\pagebreak

\section*{Appendix B}
{\footnotesize For comparison with the table in Appendix A,  we list  fees for other postgraduate taught courses in mathematics or programmes with a very strong mathematical content (what we refer to as `nearest-neighbour' programmes in the main text). 
Some of the  departments listed in Appendix A are not listed here as they do not offer such programmes.  Some  universities not listed in Appendix A  are also to be found below. For universities that are found in both tables, and for convenient cross-referencing, we have transferred the figures from the first table {\color{blue}(in blue)} into the second table and included the relative multiplicative factor {\color{red}[in red]}.  As with the previous table, this data has been collected from on-line, publicly available  information provided by  the respective institutions. Fees indicated are for the academic year 2014-2015 or 2013-2014 where information for the former is not yet available. Figures, expressed in k\pounds, have been rounded to one decimal place and courses are listed in location alphabetical order.}

\[
\hspace{-0.6cm}\begin{footnotesize}
\begin{array}{| l | l | l |l|}
\hline
 \text{
University} & \text{Fee (k\pounds) home/EU} & \text{Fee (k\pounds) non-EU} & \text{M.Sc. Course name}\\
\hline
\hline
 \text{Aberdeen}  & 3.4 & 12.6 & \text{Math.}
\\

\hline
 \text{Bath}  & 6 & 18.6 & \text{Math. Sci.}
\\
\hline
 \text{Bath}  & 6 & 18.6 & \text{Modern Appl.  Math.}
\\
\hline
 \text{Bath}  & 6 & 18.6 & \text{Math. Biol.}
\\
\hline
 \text{Birkbeck}   & 4 {\color{blue}\mbox{ }  (19)} {\color{red}\mbox{ }  [4.5]}& 7.4 {\color{blue}\mbox{ }  (20)}  {\color{red}\mbox{ }  [2.70]}& \text{Math.}
\\
\hline
 \text{Birmingham}  & 5.9 {\color{blue}\mbox{ }  (15.5)} {\color{red}\mbox{ }  [2.63]}& 13.7 {\color{blue}\mbox{ }  (15.5)} {\color{red}\mbox{ }  [1.13]}& \text{Math. Model./MORSE}
\\
\hline
 \text{Bristol}   & 8& 17.5 & \text{Math. Sci.}
\\
\hline
 \text{Cambridge (univ. + coll.)}   & 10.7 &  23.2 & \text{Part III of the Math. Tripos}
\\
\hline
 \text{Dundee}   & 3.8& 13 & \text{Math. Biol.}
\\
\hline
 \text{Durham}   & 5.7& 14 & \text{Math. Sci.}
\\
\hline
 \text{Edinburgh}    &9.3  {\color{blue}\mbox{ }  (19)}  {\color{red}\mbox{ }  [2.04]}  &17.4 {\color{blue}\mbox{ }  (22)} {\color{red}\mbox{ }  [1.26]}& \text{OR}
\\
\hline
 \text{Exeter}    & 7.5  {\color{blue}\mbox{ }  (11.9)} {\color{red}\mbox{ }  [1.59]}&17.5{\color{blue}\mbox{ }  (19.5)}  {\color{red}\mbox{ }  [1.11]}& \text{Adv. Math.}\\
\hline
 \text{Glasgow}   & 5.4 {\color{blue}\mbox{ }  (8.3)} {\color{red}\mbox{ }  [1.54]}& 17.3  {\color{blue}\mbox{ }  (17.3)} {\color{red}\mbox{ }  [1]}& \text{Math./Appl. Math.}\\

\hline
 \text{Heriot--Watt}  & 5.9 {\color{blue}\mbox{ }  (19)} {\color{red}\mbox{ }  [3.22]}& 13 {\color{blue}\mbox{ }  (22)} {\color{red}\mbox{ }  [1.69]}& \text{Math.}
\\
\hline
 \text{Heriot--Watt}  & 5.9 {\color{blue}\mbox{ }  (19)}{\color{red}\mbox{ }  [3.22]}& 13 {\color{blue}\mbox{ }  (22)}  {\color{red}\mbox{ }  [1.69]}& \text{Appl. Math. Sci.}
\\
\hline
 \text{Heriot--Watt}  & 5.9 {\color{blue}\mbox{ }  (19)}{\color{red}\mbox{ }  [3.22]}& 13  {\color{blue}\mbox{ }  (22)}  {\color{red}\mbox{ }  [1.69]}& \text{Comp. Math.}
\\
\hline
 \text{Imperial}  & 8.3 {\color{blue}\mbox{ }  (27)}  {\color{red}\mbox{ }  [3.25]}& 23 {\color{blue}\mbox{ }  (27)} {\color{red}\mbox{ }  [1.17]}& \text{Pure Math./Appl. Math.}
\\
\hline
 \text{King's}   & 8.3 {\color{blue}\mbox{ }  (22.7)} {\color{red}\mbox{ }  [2.73]}&  16.5 {\color{blue}\mbox{ }  (22.7)} {\color{red}\mbox{ }  [1.37]}& \text{Math.}
\\
 \hline
\text{Leeds}   & 5.1 {\color{blue}\mbox{ }  (9)} {\color{red}\mbox{ }  [1.76]}& 13.3 {\color{blue}\mbox{ }  (19)}  {\color{red}\mbox{ }  [1.43]}& \text{Math.}
\\

\hline
 \text{Leicester}   &8.5 {\color{blue}\mbox{ }  (8.5)} {\color{red}\mbox{ }  [1]}& 	16 {\color{blue}\mbox{ }  (16)}  {\color{red}\mbox{ }  [1]}& \text{Math. Model. Biol.}
\\
\hline
 \text{Leicester}   &8.5  {\color{blue}\mbox{ }  (8.5)} {\color{red}\mbox{ }  [1]}& 	16 {\color{blue}\mbox{ }  (16)} {\color{red}\mbox{ }  [1]}& \text{Appl. Comp.  Num. Model.}
\\
\hline
 \text{Liverpool}   & 5.3 {\color{blue}\mbox{ }  (15.4)} {\color{red}\mbox{ }  [2.90]}& 12.2 {\color{blue}\mbox{ }  (15.4)} {\color{red}\mbox{ }  [1.26]}& \text{Math. Sci.}
\\
\hline
 \text{Loughborough}   & 6.3  {\color{blue}\mbox{ }  (6.3)} {\color{red}\mbox{ }  [1]}& 13.8 {\color{blue}\mbox{ }  (13.8)}  {\color{red}\mbox{ }  [1]}& \text{Indust. Math. Model.}
\\
\hline
 \text{LSE}  & 11.6 {\color{blue}\mbox{ }  (23)} {\color{red}\mbox{ }  [1.98]}& 17.9 {\color{blue}\mbox{ }  (23)} {\color{red}\mbox{ }  [1.28]} & \text{Appl. Math.}
\\
\hline
 \text{Manchester}   & 8.4 {\color{blue}\mbox{ }  (10.7, 21.4)}  {\color{red}\mbox{ }  [1.27, 2.55]}& 14 {\color{blue}\mbox{ }  (18.4, 21.4)} {\color{red}\mbox{ }  [1.31, 1.53]}& \text{Appl. Math.}
\\
\hline
 \text{Manchester}   & 8.4 {\color{blue}\mbox{ }  (10.7, 21.4)} {\color{red}\mbox{ }  [1.27, 2.55]}& 14 {\color{blue}\mbox{ }  (18.4, 21.4)} {\color{red}\mbox{ }  [1.31, 1.53]}& \text{Pure Math.  Math. Logic}
\\

\hline
 \text{Nottingham}   & 6  {\color{blue}\mbox{ }  (6)} {\color{red}\mbox{ }  [1]}& 13.7 {\color{blue}\mbox{ }  (13.7)} {\color{red}\mbox{ }  [1]}& \text{Pure Math.}
\\
\hline
  \text{Oxford}    & 5.6 {\color{blue}\mbox{ }  (24.9, 27.5)}  {\color{red}\mbox{ }  [4.45, 4.91]}& 15.3 {\color{blue}\mbox{ }  (24.9, 31.1)}   {\color{red}\mbox{ }  [1.62, 2.03]}& \text{Math. Model. Sci. Comp.}
\\
\hline
 \text{Queen Mary}   &6.3 {\color{blue}\mbox{ }  (15.5)} {\color{red}\mbox{ }  [2.46]}&  13.5 {\color{blue}\mbox{ }  (18)} {\color{red}\mbox{ }  [1.33]}& \text{Math.}
\\

\hline
 \text{Sheffield}   &  6 {\color{blue}\mbox{ }  (6)} {\color{red}\mbox{ }  [1]}&15.4 {\color{blue}\mbox{ }  (15.4)} {\color{red}\mbox{ }  [1]}& \text{Math.}
\\
\hline
 \text{Southampton}   &  7.3 & 15 & \text{OR}
\\
\hline
 \text{Surrey}   &  6.3 & 16.3 & \text{Math.}
\\
\hline
 \text{Sussex}   &  5.5 & 13 & \text{Math.}
\\
\hline
 \text{UCL}   &8.5 {\color{blue}\mbox{ }  (19.6)} {\color{red}\mbox{ }  [2.30]}&16.8 {\color{blue}\mbox{ }  (21.7)} {\color{red}\mbox{ }  [1.29]}& \text{Math. Model.}
\\
\hline
\text{Warwick}   & 7{\color{blue}\mbox{ }  (31)} {\color{red}\mbox{ }  [4.42]}& 15.9  {\color{blue}\mbox{ }  (31)} {\color{red}\mbox{ }  [1.94]}& \text{Math.}\\
\hline
 \text{York}   & 6.2 {\color{blue}\mbox{ }  (17.4, 16.5)} {\color{red}\mbox{ }  [2.80, 2.66]}& 14.3 {\color{blue}\mbox{ }  (22.2, 21)}  {\color{red}\mbox{ }  [1.55, 1.47]}& \text{Adv. Math. Biol.}
\\

%\hline
%8 &\text{LSE}   & 23 & 23 & \text{Statistics (Financial Statistics)}
%\\
%\hline
%\text{Glasgow}  & 8.3& 17.2 & \text{Financial Modelling}\\

\hline

\end{array}
\end{footnotesize}\]

\newpage
\noindent {\footnotesize Based on the figures in red above, we group universities according to the differential between the presented M.Sc. in Financial Mathematics and other M.Sc. programmes in a strongly mathematical subject area on the same campus.}

\[
\begin{footnotesize}
\begin{array}{| l | c | c| c| c| c| }
\hline
\text{Differential} & 0-10\% & 10-20\% & 20-50 \% & 50-100\% & > 100\% \\
\hline
\hline
\text{Home-EU} 
& 
\begin{array}{l}
\text{Leicester}\\
\text{Loughborough} \\
 \text{Nottingham}\\
\text{Sheffield} 
\end{array}
& 
& 
\begin{array}{l}
\text{Glasgow} \\
\text{Manchester}
\end{array}
& 
\begin{array}{l}
\text{Exeter} \\
\text{Leeds} \\
\text{LSE}
\end{array}
&
\begin{array}{l}
\text{ }\\
\text{Birkbeck}\\
\text{Birmingham}\\
%\text{Durham}\\
\text{Heriot--Watt}\\
\text{Imperial}\\
\text{King's}\\
\text{Liverpool}\\
\text{Manchester}\\
\text{Queen Mary}\\
\text{UCL}\\
\text{Warwick}\\
\text{York}\\
\text{ }
\end{array}  \\
\hline
\text{Overseas} & 
\begin{array}{l}
\text{Glasgow}\\
\text{Leicester}\\
\text{Loughborough}\\
\text{Nottingham}\\
\text{Sheffield}
\end{array}
&
\begin{array}{l} 
\text{Birmingham}\\
\text{Exeter}\\
\text{Imperial}
\end{array}
& 
\begin{array}{l}
\text{ }\\
\text{Edinburgh}\\
\text{King's}\\
\text{Leeds}\\
\text{Liverpool}\\
\text{LSE}\\
\text{Manchester}\\
\text{Queen Mary}\\
\text{UCL}\\
\text{York}\\
\text{ }
\end{array}
& 
\begin{array}{l}
\text{Heriot-Watt}\\
\text{Oxford}\\
\text{Warwick}\\
\text{York}
\end{array}
&
\begin{array}{l}
\text{Birkbeck}\\
\text{Oxford}
\end{array}  

\\
\hline

\end{array}
\end{footnotesize}
\]

\end{document}